\documentclass[12pt]{article}
\usepackage{amsmath,amssymb,amsthm,xspace,a4wide,amscd}
\newtheorem{theorem}{Theorem}%[section]
\newtheorem{lemma}{Lemma}%[section]
\newtheorem{proposition}{Proposition}%[section]
\newtheorem{corollary}{Corollary}%[section]
 \newtheorem{remark}{Remark}

\renewcommand{\phi}{\varphi}
\begin{document}
\title{Imaginary Verma modules for the  extended affine Lie algebra $sl_{2}(\mathbb C_{q})$}
\author{M. Dokuchaev, L. M. Vasconcellos Figueiredo \\
 and V. Futorny \footnote{Regular Associate of the ICTP}}

\date{}
\maketitle

\begin{abstract}
We consider one of the most natural extended affine Lie algebras , the
 algebra
 $sl_2({\mathbb C}_q)$, over the quantum torus 
${\mathbb C}_q$ and begin a theory of its representations. In particular, we study a class of imaginary
Verma modules, obtain for them a criterion of irreducibility and describe their submodule structure when 
the modules are in ``general position''.    
\end{abstract}

\centerline{Mathematical Subject Classification 17 B67, 17 B65}

\newpage

\section { Introduction}
\hspace{0.4cm}
\medskip

In  recent years Extended affine Lie algebras (EALA's)
 have been studied extensively (see \cite{AABGP}
and references therein). These algebras provide a natural generalization of affine Kac-Moody Lie
algebras and more generally of toroidal Lie algebras. Certain representations 
of toroidal  
Lie algebras were considered in \cite{BC}, \cite{CF}, \cite{FK}. In particular, {\it imaginary 
Verma modules} were studied in \cite{FK}. For affine Lie algebras these modules were studied 
in \cite{F} and in \cite{C}. Some vertex operator representations for toroidal Lie algebras 
were constructed in \cite{EMY1}, \cite{EMY2}, \cite{EM}, 
\cite{B}, \cite{BB}, \cite{T1}, \cite{T2}. Also in \cite{BB} certain
 induced irreducible  modules were considered.  

Vertex operator representations for  extended affine Lie algebras were 
discussed in \cite{BS}, \cite{G1}, \cite{G2}, \cite{T3}, \cite{W}, \cite{Y}.
A class of extended affine Lie algebras coordinatized by a quantum torus (cf. [M])
 contains in particular 
all toroidal and all affine Lie algebras.  Such algebras were considered in 
 \cite{BS}, \cite{G1} and \cite{G2}.    
 
In the present paper we consider an extended affine Lie algebra $sl_2({\mathbb C}_q)$ over the 
{\it quantum torus} ${\mathbb C}_q$ and study  imaginary Verma modules for this algebra. 
Our results generalize the results in \cite{FK} about the imaginary Verma modules. In particular, 
we obtain a criterion of irreducibility for such modules and  describe the submodule structure and irreducible 
quotients in the
case when at least one central element acts non-trivially. 

The structure of the paper is the following. In Section 2, we provide all necessary preliminaries 
on the quantum torus and our algebra. In Section 3, we study  Verma modules for the 
generalized Heisenberg subalgebra.  Theorem 1 describes the irreducible quotient of a Verma module for 
the generalized Heisenberg subalgebra when a particular central element acts non-trivially. 
Finally, in Section 4, we study  imaginary Verma modules.   Here, Corollary 1 establishes a criterion 
of irreducibility for imaginary Verma modules. Theorem 2 describes the submodule structure of imaginary Verma modules 
in the case when at least one central element acts non-trivially. Corollaries 2 and 3 give a description 
of the irreducible quotients in this case.

\bigskip
\section{Preliminaries} 

 \subsection{The Quantum Torus}

  Let $q \in M_{n}(\mathbb C)$ be a $n\times n$ matrix $q=(q_{ij})$, 
$1\leq i,j \leq n$, such that 
$q_{ii}=1$ with $q_{ij}=
q_{ji}^{-1}$. Let $\mathbb C_{q}$ be the complex algebra spanned by $t_{i}^{\pm 1}$, $1\leq i \leq n$ and the
product $t_{i}t_{j}$ defined by $t_{i}t_{j}=q_{ij}t_{j}t_{i}$ for all $i,j$ with $1\leq i,j \leq n$.
 Note that for all $k,l \in 
\mathbb Z$ , 
$$t_{i}^{k}t_{j}^{l}=q_{ij}^{kl}t_{j}^{l}t_{i}^{k}.$$

 Let $\Lambda = \mathbb Z^{n}$ and $a\in \Lambda$ , $a=(a_{1}, a_2,...,a{_n})$ and 
define $t^{a}= t_{1}^{a_{1}}t_{2}^{a_{2}}...t_{n}^{a_{n}}$ and $\mathbb C^{a}_{q}=
\mathbb C t^{a}$. For $a,b \in \Lambda$ we have that
$$ t^{a}t^{b} =\sigma(a,b)t^{a+b},$$

\noindent where  $ \sigma:\Lambda \times \Lambda \rightarrow \mathbb C $ is defined by 
    $$\sigma(a,b)= \prod_{1\leq i<j \leq n}q_{ji}^{a_{j}b_{i}}.$$
   Then  $$\sigma(b,a)= \prod_{1\leq i<j \leq n}q_{ji}^{b_{j}a_{i}}.$$

and 
$$\sigma(b,a)^{-1}= \prod_{1\leq i<j \leq n}q_{ij}^{a_{i}b_{j}}.$$
Here we use the notation from \cite{BGK}. Define 
$$f(a,b)=\sigma (a,b)\sigma (b,a)^{-1}=\prod_{i\neq j}q_{ij}^{a_{i}b_{j}}.$$

 \noindent     Note that $[t^{a},t^{b}] = t^{a}t^{b}-t^{b}t^{a} = \sigma(a,b)t^{a+b}-
\sigma(b,a)t^{a+b} = \sigma(b,a)(f(a,b)-1)t^{a+b}.$

     The basic properties of $f$ are:
\begin{proposition}\label{p1} For all $a,b,a',b' \in \Lambda \,$ we have:
\begin{enumerate}
\item  $f(a+a',b)=f(a,b)f(a',b);$
\item $f(a,b+b')=f(a,b)f(a,b');$
\item $f(b,a)= f(a,b)^{-1};$
\item $f(a,a)=1$; $f(a,-a)=1;$ 
\end{enumerate}

\end{proposition}

  Define the {\it radical }of $\mathbb C_{q}$ by

$$ R= \{ a\in \Lambda \, |\,  f(a,b)=1\,  \forall  b\in \Lambda \}.$$

{\bf  Remarks:}
\begin{enumerate}
\item  Note that the algebra $\mathbb C_{q}$ is  commutative if and only if 
$q_{ij}=1$ for all $i,j$ with $1\leq i,j \leq n$, or equivalently if $R =\Lambda.$

\item  The {\it center} of $\Bbb C_{q} $ is
 $$ Z(\mathbb C_{q}) = \sum_{a\in R}\Bbb C t^{a} $$ and
 $$ \mathbb C_{q} = Z(\mathbb C_{q})\oplus [\mathbb C_{q},\mathbb C_{q}].$$ 
 \noindent To see this just observe that $t^{a}\in [\mathbb C_{q},\mathbb C_{q}]$ if and only
if $a \notin  R$.

\end{enumerate}

     Now let $${\mathfrak g} = sl_{2}(\Bbb C_{q}) 
=\left \{  \left( \begin{array}{cc}
x & y \\ 
 z & w \\ \end{array} \right) \,|\, x,y,z,w \in \Bbb C_{q}, \,x+w \in [\Bbb C_{q},\Bbb
C_{q}] \right \}.$$
      
      A basis for $sl_{2}(\Bbb C_{q}) $ is the set 
$$\left \{ x_{a} ,y_{a} , u_{a}, w_{b}\,|\,a,b \in \Lambda , b \notin R \right \}$$

where  $$x_{a}=\left( \begin{array}{cc} 0 & t^{a}\\ 0 & 0 \end{array} \right),
 y_{a} =\left(\begin{array}{cc}  0&  0 \\ t^{a} & 0 \end{array}\right),
u_{a} = \left( \begin{array}{cc} t^{a}  & 0\\ 0 & -t^{a}\end{array}\right),
w_{b}=\left( \begin{array}{cc} t^{b}  & 0\\ 0 & t^{b} \end{array}\right) .$$ 
      We have that $\mathfrak g$ is a  $\Bbb C $ -Lie algebra with brackets $[\, ,\,]_{old}$ being defined by $[x,y]_{old} = xy-yx $ for $x,y \in  \mathfrak g$.

      Call, for $a \in \Lambda$,

 $${\mathfrak{g}^{a}}  =\left\{ \begin{array}{ll} {\mbox{span}}\{ x_{a}, y_{a}, u_{a} , w_{a} \} & {\mbox{if}} \; a \notin R, \\  {\mbox{span}}\{ x_{a}, y_{a}, u_{a}  \} & {\mbox{if}} \; a \in R.\\ \end{array} \right.$$

%\right

       For $ i=1,2,..,n,$  let $d_{i} $ be the {\it derivations} on 
$\mathfrak g $   defined by  
  $$ d_{i}x = a_{i}x , x \in\frak g^{a}, a = (a_{1}, a_{2},..., a_{n})\in \Lambda.$$

 Define  $\epsilon : \Bbb C_{q}\rightarrow \Bbb C$ by

$$\epsilon(t^{a}) =\left\{\begin{array}{ll} 0 & \mbox{if } \, a\neq (0,0,...,0),\\  1 & \mbox{if } \, a= (0,0,...,0) \end{array} \right.$$

\noindent and define an invariant  non-degenerate bilinear form $( \,,\, )$ on $\mathfrak g$  by $(x,x')=\epsilon(tr(xx'))$ where $tr(\, \,)$
indicates the usual trace function.

Let $C =\bigoplus_{i=1}^{n}\Bbb C c_{i}$ , $D=\bigoplus_{i=1}^{n} \Bbb C d_{i}$ and let $\mathfrak q=\mathfrak g \oplus C\oplus D$ be the Lie algebra  with brackets

$$[x,x']:= [x,x']_{new} = [x,x']_{old} +\sum _{i=1}^{n} ( [d_{i},x],x')c_{i}$$ for
$x,x'\in \frak g$. 

$$ [c_{i}, x]=0 \forall x\in \mathfrak q ,\,\forall  i=1,2,...,n$$  and
$$ [D,D]=0 , [d_{i},x] = d_{i}x, \,  \forall x\in \frak g \, ,\forall i=1,2,...,n .$$
Then $\mathfrak q$ is a Lie algebra, and it is the {\it quantum torus }.

The brackets of the basis elements of $\mathfrak q $ are:
$$[x_{a},y_{b}]=\frac{\sigma(b,a)(f(a,b)+1)}{2}u_{a+b}+\frac{\sigma(b,a)(f(a,b)-1)}{2}w_{a+b}+\delta_{a,-b}\sigma(a,b)\sum_{i=1}^{n}a_{i}c_{i} $$
$\forall a,b \in \Lambda$, and where $\delta_{a,-b} = 0$ if $a\neq -b$ and $\delta_{a,-a} = 1$, 
$$[x_{a}, w_{b}]=\sigma(b,a)(f(a,b)-1)x_{a+b}\,\,  \forall a,b \in \Lambda ,b\notin R$$
$$[x_a,u_b]=-\sigma(b,a)(f(a,b)+1)x_{a+b} \,\,\forall a,b\in \Lambda$$
$$[y_a,w_b] =\sigma(b,a)(f(a,b)-1)y_{a+b} \,\,\forall a,b\in \Lambda ,b \notin R$$

$$[y_a,u_b]=\sigma(b,a)(f(a,b)+1)y_{a+b}\,\, \forall a,b\in \Lambda$$
$$[w_a, w_b]=
[u_{a},u_{b}]=\left \{ \begin{array}{ll} 2\sigma(a,-a)\sum_{i=1}^{n}a_{i}c_{i} & \mbox {if $b=-a$}\\ 0&\mbox{if $a+b \in R ,a+b \neq 0$}\\   \sigma(b,a)(f(a,b)-1)w_{a+b} & \mbox{if $a+b \notin R$} \end{array}\right.$$
$$[u_{a},w_b]=\left \{ \begin{array}{ll} 0&\mbox{if $a+b\in R$}\\ \sigma(b,a)(f(a,b)-1)u_{a+b}& \mbox{if $a+b \notin R.$} \end{array} \right. $$

\subsection{Roots}
 
 Let $\mathfrak h = \Bbb C h\oplus C \oplus D$ where  $ h = \left( \begin{array}{rr} 1&0\\ 0& -1 \end{array} \right)$. Then $ \mathfrak h$
is a Cartan Subalgebra of $\mathfrak q$ . 
  
We shall describe the root space decomposition of $\mathfrak q$ with respect
 to $\mathfrak h$.

Let $\alpha \in \mathfrak h ^{*}$ be such that $\alpha(h)=2$ and $\alpha|_{C,D}
=0$.  For all $i,j=1,2,...,n$ let $\delta_i \in \mathfrak h^{*}$ be defined by
$\delta_i(c_j)=0$,  $\delta_i(h)=0$,  $\delta_i(d_j)=\delta_{i,j}$ for all $ j=
1,2,...,n$.
 We identify $a=(a_{1}, a_{2},..., a_{n}) \in \Lambda$  with the element 
$\sum _{i=1}^{n} a_{i}\delta_{i} \in \mathfrak h^{*}$ so that for all $x\in \mathfrak g^a $ 
$$[d_i,x]=a_i x = (\sum _{i=1}^{n}a_i \delta_i)(d_i)x =a(d_i)x $$ for all 
$i=1,2,....n .$
Then the {\it root system} $\Delta$ of $\mathfrak q$ is
$$ \Delta = \{ \pm \alpha + \Lambda \} \cup  \Lambda \backslash \{0\}.$$
The {\it root space decomposition} of $\mathfrak q$ with respect to $\mathfrak h$  is  $$ \mathfrak q = \mathfrak h \oplus \sum _{\beta \in \Delta} \mathfrak q_{\beta} = \sum _{\beta \in \Delta \cup \{ 0\}} \mathfrak q_{\beta},$$
where $$\mathfrak q_{\beta}= \Bbb C x_a $$ if $\beta =\alpha +a ,a\in \Lambda$
$$\mathfrak q_{\beta}= \Bbb C y_a $$ if $\beta =-\alpha +a$ if $ a\in \Lambda$

$$\mathfrak q_{\beta}= \Bbb C u_a\oplus \Bbb C w_a $$ if $\beta=a  \, ,a\notin R$ 
 
$$\mathfrak q_{\beta}= \Bbb C u_a$$ if $\beta =a, a\in R$.

\section {Heisenberg subalgebra}

  Let  $ \mathfrak t$ be the {\it generalized Heisenberg subalgebra } of  $\mathfrak q$ , 
$$ \mathfrak t = \sum_{a \in \Lambda \setminus \{ 0\}} \mathfrak q_a \oplus C .$$

 We can order the elements of $\Lambda$ lexicographically, that is,
 for $a,b \in \Lambda, a= (a_1, a_2,..., a_n)$ and $b= (b_1,b_2,...,b_n)$, 
$a < b$ if and only if, for some $ i=1,...,n$, $a_i < b_i$, and  for all  $j > i$,
 $a_j=b_j$. Set $\Lambda^- =\{a\in \Lambda | a < 0=(0, \ldots, 0)\}$.

Consider the following Lie subalgebras of $\mathfrak t$,
$$ \mathfrak t^{\pm} = \sum_{\stackrel{ a>0}{(a<0)}} \mathfrak q_a $$ and $$\mathfrak b =
C \oplus \mathfrak t^{+}, $$ a {\it Borel subalgebra} of $\mathfrak t$.
 
 Fix $\lambda \in C^{*}$. Let $H(\lambda)$ be the {\it Verma-type 
module  for $\mathfrak t$  associated to} $\lambda $, that is 
$$ H(\lambda) = {U(\mathfrak t)}/{ <\mathfrak t^{+}, c_i -\lambda(c_i)1,i=1,...,n >},$$
\noindent where $U(\,\,)$ denotes the universal enveloping algebra. From the
general theory of Verma modules for Lie algebras we know that if $V = \Bbb C v_{\lambda}$
is a one-dimensional vector space and if we make $V$ a $\mathfrak b$- module
as $\mathfrak t^{+} v_{\lambda} =0$  and for $ i=1,2,...,n$, $c_{i} v_{\lambda}= \lambda (c_{i})v_{\lambda}$, \,then
$$H(\lambda) \cong  U(\mathfrak t) \otimes _{U(\mathfrak b)}\Bbb C v_{\lambda}.$$  Because 
$U(\mathfrak t) =U(\mathfrak t^{-}) \otimes U(C) \otimes U(\mathfrak t^{+})$ we
have that $$ H(\lambda) \cong U(\mathfrak t^{-}) \otimes _ { \Bbb C }\Bbb C v_{\lambda}$$
\noindent as vector spaces.

We remark that the elements $d_{i}, i= 1,...,n$ of $\mathfrak q$ act on $H(\lambda)$ as 
$$d_{i}.v_{\lambda} = 0$$
\noindent and if $r\geq 1$ and if for $j=1,...r , z_{a_{j}} \in \{u_{a_{j}}, w_{a_{j}} \},$

$$d_{i}.(z_{a_{1}}...z_{a_{r}}v_{\lambda}) = \sum_{j=1}^r z_{a_{1}}z_{a_{2}} ... [d_{i},z_{a_j}] ... z_{a_{r}} v_{\lambda},$$\noindent defining a $\Lambda$-gradation on $H(\lambda)$
$$H(\lambda) = \sum _{\beta \in \Lambda}H(\lambda)_{\beta},$$
\noindent where, if $\beta=(\beta_{1},...,\beta_{n}),$

$$H(\lambda)_{\beta} = \{x \in H(\lambda) | d_{i}.x = \beta_{i}x \,\forall i=1,...,n\}.$$

\noindent If $x\in H(\lambda)_{\beta}$ we say that $x$ is {\it homogeneous of degree }$\beta \in \Lambda.$
      
     Our goal is to describe the irreducible quotients of $H(\lambda)$.
     For $a = (a_{1},a_{2},...,a_{n}) \in \Lambda$ define 
$$\mathfrak t_{(a)} = \sum_{l\in \mathbb Z \setminus \{0 \}} \mathfrak q_{la} \oplus 
\mathbb C (a_{1}c_{1} + a_{2}c_{2} + ... + a_{n}c_{n}).$$    
\noindent Then $\mathfrak t_{(a)}$ is a Heisenberg subalgebra of $\mathfrak t$ as in the affine case. Consider the 
$\mathfrak t_{(a)}$-submodule $H_{(a)}$ of $H(\lambda)$ generated by $v_{\lambda}$. Form 
    $$ \hat H_{(a)} = \sum_{\beta \in \Lambda \setminus \{ (0,...,0) \} }(H_{(a)})_{\beta}.$$
 \noindent Let $\lambda \in C^{*}$ and let
   $$\Lambda_{\lambda} =\{ (r_{1},r_{2},...,r_{n}) \in \Lambda \, | \, r_{1}\lambda (c_{1}) +r_{2} \lambda (c_{2})+...+r_{n} \lambda (c_{n}) = 0,
(r_1,r_2,...,r_n) \neq (0,0,...,0) \}.$$ 
    For each $r$ with $1 \leq r \leq n $ define
    $$R_{r} = \{a = (a_{1},...,a_{r},0,...,0) \in \Lambda \,| \, f(a,b) = 1\, \mbox{for all} \, b = (b_{1},...,b_{r},0,...,0) \in \Lambda \}.$$

\begin{proposition}\label{p2} Let  $a \in \Lambda_{\lambda} \cap R_{r} $. The $\mathfrak t$-submodule $\tilde H_{(a)} = U(\mathfrak t)\hat H_{(a)}$ is a proper
submodule of $H(\lambda)$.
 \end{proposition}

\begin{proof} Assume in  contradiction that $\tilde H_{(a)} = H(\lambda)$. Then $v_{\lambda} \in \tilde H_{(a)}$ and
there exist a monomial $x \in \hat H_{(a)}$ and a monomial $y \in U(\mathfrak t)$ such that $v_{\lambda}$ appears in
 the support of $yx$. In particular, $yx \neq 0$. We can write $x = x_{l_{1}a}x_{l_{2}a}...x_{l_{m}a}v_{\lambda}$
where each $l_{i}a < 0$ and $x_{l_ia} \in \{u_{l_ia},w_{l_ia}\}$. Then $x \in  (\hat
H_{(a)})_{\beta}$ with $\beta = \sum_{i=1}^{m}l_{i}a$. Let $y$ be of degree
$\gamma$. If $\gamma < -\beta$ then $\gamma +\beta < 0$ and $yx \in (\tilde{H}_{(a)})_{\gamma +\beta} \not\ni v_{\lambda}$,
which is not possible. Hence $\gamma \geq -\beta > 0$.
Write, without loss of generality, $y = z_{b}z_{b'}...z_{\tilde b}$ where $b \leq b' \leq ...\leq \tilde b$  and $z \in \{ u, w \}$. Then $\gamma = b + b' + ...+ \tilde b$ and clearly, $\tilde b > 0$. If $\tilde b = (\tilde b_{1},...,\tilde b_{s}, 0,...,0)$
with  $\tilde b_{s} \neq 0$ and $s > r $ then $ \tilde b + l_{i}a > 0 $ for all $i$ and $z_{\tilde b}x = 0$ which implies that $yx = 0 .$ So $s \leq r.$ 

If $\tilde b$ is different from each $-l_{i}a (i=1,...,m)$, then $[z_{\tilde b},x_{l_{i}a}] = 0$ as $-l_{i}a \in R_{r} , (i = 1,...,m)$. It follows that $z_{\tilde b}x = z_{\tilde b}x_{l_{1}a}...x_{l_{m}a}v_{\lambda} = x_{l_{1}a}...x_{l_{m}a}z_{\tilde b}v_{\lambda} = 0 $, and thus $yx =0 $, a contradiction. Hence 
we may assume that $\tilde b = -l_{i}a$ for a fixed $i$. Then $$[z_{\tilde b}, x_{l_{i}}]v_{\lambda} = [z_{-l_{i}a}, x_{l_{i}a}]v_{\lambda}
= -2 \sigma (-l_{i}a, l_{i}a)l_{i}\sum _{j=1}^{n}a_{j}\lambda (c_{j})v_{\lambda} =0 $$ as $l_{i}a  \in \Lambda_{\lambda}.$ Moreover $[z_{\tilde b}, x_{l_{j}a} ] = 0$ for all $j$ with $l_{j} \neq l_{i}$, because $l_{i}a \in R_{r}$. Then again we obtain that $z_{\tilde b} x = 0$, the
 final contradiction.

\end{proof}

Set $$ \tilde H = \sum _{a \in \Lambda_{\lambda} \cap (\cup_{r=1}^{n} R_{r})} \tilde H_{(a)}. $$

\noindent It follows from Proposition~\ref{p2}  that $\tilde H$ is a proper  $\mathfrak t$ - submodule of $H(\lambda)$.                                                                                                                                                                                                                                                                                                                                                                                                                                                                                                                                                                                                                                                                                                                                                                                                                                                                                                                                                                                                                                                                                                                                                                                                                                                                   
Next is the main result in this section.                

\begin{theorem} \label{t1} Let $\lambda \in C^{*}.$
     If $\lambda (c_{1}) \neq 0 $ then $H(\lambda)/{\tilde H}$ is irreducible.

\end{theorem}

   In order to prove the theorem we shall need several lemmas. 

\begin{lemma}\label{l1} Define an ordering $\preceq$ on $\Lambda^{-}$ by:

\begin{enumerate}

\item $(-k,0,...,0) \prec (0,...,0)$ if $k > 0;$  $(-k,0,...,0) \prec (-m,0,...,0)$ with $k,m > 0$ if $k < m.$

For $r > 1,$

\item $(k_{1},...k_{r-1},-k,0,...,0) \prec (m_{1},...,m_{r-1},0,...,0)$ if $k > 0.$
\item $(k_{1},...,k_{r-1}, -k,0,...,0) \prec (m_{1},...,m_{r-1},-m,0,...,0)$ with $k,m > 0$ if either $k < m $ or $k = m$ and $(k_{1},..., k_{r-1},0,...,0) < (m_{1},..., m_{r-1},0,...,0)$ in the previous ordering.
\end{enumerate}
 Define an ordering on the basis $ \{ w_{a}, u_{a} , a\in \Lambda ^{-} \}$ of $\mathfrak t^{-}$ putting $w_{a} < u_{a}$ for $a\in \Lambda^{-}$
and $z_{a} < z_{b}$ if $a \prec b$ , where $z_{a} \in \{ u_{a}, w_{a} \}, z_{b} \in \{u_{b}, w_{b} \}.$

\noindent The monomials $$z_{a_{1}}...z_{a_{s}}v_{\lambda}$$ with $a_{i} \in \Lambda^{-},$  and $z_{a_{1}}\leq z_{ a_{2}}\leq ... \leq z_{a_{s}}$ form a basis of $H(\lambda)$.      
\begin{proof} Just apply the P-B-W Theorem.
\end{proof}
\end{lemma}

\begin{lemma}\label{l2} Let $b = (b_{1},...,b_{r},0,...,0) \in \Lambda \setminus R_{r}$. Then for arbitrary positive integers $k_{1},...,k_{r-1}$ there exist integers
$N_{i} > k_{i}$ $(i=1,...,r-1)$ such that $$f((-N_{1},-N_{2},...,-N_{r-1},1,0,..,0), (b_{1},b_{2},...,b_{r},0,...,0)) \neq 1.$$
\begin{proof}Observe first that an element $x = (x_{1},...,x_{r},0,...,0) \in \Lambda$ is in  $R_{r}$ if and only if

\[ (1) \, \,\,  \left \{  \begin{array}{lllll}

         q_{21}^{x_{2}}q_{31}^{x_{3}}q_{41}^{x_{4}}...q_{r1}^{x_{r}} & =  1  \\

q_{21}^{-x_{1}}q_{32}^{x_{3}}q_{42}^{x_{4}}...q_{r2}^{x_{r}} & = 1 \\

 q_{31}^{-x_{1}}q_{32}^{-x_{2}}q_{43}^{x_{4}}...q_{r3}^{x_{r}} & =  1 \\      
  .........................\nonumber\\

q_{r1}^{-x_{1}}q_{r2}^{-x_{2}}q_{r3}^{-x_{3}}...q_{r r-1}^{-x_{r-1}} & =  1.   \\

\end{array} \right. \]

\noindent Indeed, $x \in R_{r}$ if and only if $f(x,b) =1$ for all $b = (b_{1},...,b_{r},0,...,0) \in \Lambda$. It is easy to see that 
$$1 = f(x,b)$$    $$ = (q_{21}^{x_{2}}q_{31}^{x_{3}}q_{41}^{x_{4}}....q_{r1}^{x_{r}})^{b_{1}}(q_{21}^{-x_{1}}q_{32}^{x_{3}}q_{42}^{x_{4}}....q_{r2}^{x_{r}})^{b_{2}}(q_{31}^{-x_{1}}q_{32}^{-x_{2}}q_{43}^{x_{4}}....q_{r3}^{x_{r}})^{b_{3}}...(q_{r1}^{-x_{1}}q_{r2}^{-x_{2}}q_{r3}^{-x_{3}}....q_{rr-1}^{-x_{r-1}})^{b_{r}} .$$

\noindent Since this holds for all $b_{1},...,b_{r}$, it follows that $x \in R_{r}$ if and only if $x$ satisfies $(1)$. Now let $N_{1},...,N_{r-1}$ be positive integers with $N_{i} > k_{i} (i = 1,...,r-1)$. Then the equality $f((-N_{1},...,-N_{r-1},1,0,...,0), (b_{1},...,b_{r},0,...,0)) = 1$ implies
$$(q_{21}^{b_{2}}q_{31}^{b_{3}}q_{41}^{b_{4}}....q_{r1}^{b_{r}})^{N_{1}}(q_{21}^{-b_{1}}q_{32}^{b_{3}}q_{42}^{b_{4}}....q_{r2}^{b_{r}})^{N_{2}}(q_{31}^{-b_{1}}q_{32}^{-b_{2}}q_{43}^{b_{4}}....q_{r3}^{b_{r}})^{N_{3}}...(q_{r-11}^{-b_{1}}q_{r-12}^{-b_{2}}q_{r-13}^{-b_{3}}....q_{rr-1}^{b_{r}})^{N_{r-1}}$$ $$= q_{r1}^{-b_{1}}q_{r2}^{-b_{2}}q_{r3}^{-b_{3}}...q_{r r-1}^{-b_{r-1}}.$$

\noindent If this holds for all $N_{1},..., N_{r-1}$ with $N_{i} > k_{i}$ then each factor in brackets has to be 1, which gives $(1)$. So
$b\in R_{r}$, a contradiction which proves the lemma.

 \end{proof}
\end{lemma}

\begin{lemma}\label{l3} Let $0 \neq x \in H(\lambda)$ be  homogeneous  of degree $\alpha \in \Lambda,$ $ \alpha = (\tilde \alpha, -L,0,...,0),$ $ L > 0, $
$\tilde \alpha \in \mathbb Z^{r-1}.$ Suppose that  at least one monomial of $x$ contains a factor
$z_{a'}$ with $a' = (a'_{1},...,a'_{r},0,...,0) \notin R_{r}$ and $a'_{r} \neq 0$. If $x \notin \tilde H $ then there exists $y \in \mathfrak t$ such that $yx \notin \tilde H$ and $yx \in H(\lambda)_{\beta}$ with $\beta = (\tilde \beta, -L+1,0,...,0)$.

\begin{proof}  We assume that all monomials $ M = z_{a_{1}}z_{a_{2}}...z_{a_{s}}v_{\lambda}$ ($z \in \{u,w\}, a_{i}\in \Lambda^{-})$,
are ordered as in Lemma~\ref{l1}. If some $a_{i} \in R_{r}$, then $f(a_{i}, a_{j}) = 1$ for all $j=1,...,s$, so that we can 
write $M = z_{a_{i_{1}}}....z_{a_{i_{t}}}z'v_{\lambda}$ where $z'= z_{a_{j_{1}}}...z_{a_{j_{s-t}}}$ with $a_{j_{l}} \in R_{r}.$
We say that such monomial has length $t+1$ (so the length of $z'v_{\lambda}$ is 1). 
Take an arbitrary ordering on the set of all monomials $z'v_{\lambda}$ and, using it, order the set of monomials                         $z_{a_{1}}z_{a_{2}}...z_{a_{s}}z'v_{\lambda}$ of a fixed length $s+1$ lexicographically from left to right. 
Write $$x = \sum_{i\in I}\lambda_{i}z_{i}v_{\lambda}$$
\noindent with $$z_{i} = z_{a_{1i}}^{t_{1i}}z_{a_{2i}}^{t_{2i}}...z_{a_{s_{i}i}}^{t_{s_{i}i}}z'_{i},$$   where each $t_{li} \geq 1$ and 
$a_{j_{1}i} \neq a_{j_{2}i}$ if $j_{1} \neq j_{2}$.    
  Among all monomials $z_{i}v_{\lambda}$ consider those which have  maximal length and denote the corresponding subset of indices by $I'$.    
Let $z_{i_{0}}v_{\lambda}$ be the smallest monomial in $ X' = \{ z_{i}v_{\lambda}\, | \,  i\in I' \}$ and set   
$X = \{z_{i}v_{\lambda}\,|\,i\in I\}.$

For each $i \in I$ and $1 \leq l \leq s_{i}$, write $a _{li} = ( \tilde a_{li}, -m_{li},0,...,0)$, where 
$m_{li} \geq 0, \, \tilde a_{li} = $ $(k_{li}^{(1)}, k_{li}^{(2)},..., k_{li}^{(r-1)})$ and set $ k_{j} = \mbox{max}_{i\in I,1 \leq l \leq s_{i}}
\{ |k_{li}^{(j)}| \}$. 
We consider separately the cases $m_{1i_{0}} > 1 $ and $m_{1 i_{0}} = 1$. 
So suppose first that   $m_{1i_{0}} > 1 $.
By Lemma~\ref{l2}, there exist $N_{1},...,N_{r-1}$ with $N_{j} > k_{j}$ such that 
$$f((-N_{1},...,-N_{n-1},1,0,...,0),(\tilde a_{1i_{0}},-m_{1i_{0}},0,...,0)) \neq 1,$$
 so setting $\bar b = (-N_{1},...,-N_{r-1})$ and $b = (\bar b,1,0,...,0)$ we have that $f(b, a_{1i_{0}}) \neq 1.$ 
In particular,
$b+a_{1i_{0}} \notin R_{r} $ as $f(b+a_{1i_{0}}, a_{1i_{0}}) =f(b, a_{1i_{0}}).$ 
Taking $y = w_{b}$ we will show that $0 \neq yx \notin \tilde H$. 
First we calculate $yz_{i}v_{\lambda}$ for arbitrary $i \in I$. 
Set $\xi_{ji} = \sigma(a_{ji},b)(f(b,a_{ji})-1)$. 
Because of the choice of $b$, $a_{ji} \neq -b$ for all $i\in I$, $j= 1,2,...,s_{i}$. 
Hence $[w_{b}, z_{a_{ji}}] = \xi_{ji}z_{b+a_{ji}}.$ We have

$$yz_{i}v_{\lambda}= \xi_{1i}t_{1i}z_{(\tilde a_{1i}+\bar b,-m_{1i}+1)} z_{(\tilde a_{1i},-m_{1i})}^{t_{1i}-1}z_{a_{2i}}^{t_{2i}}...z_{a_{s_{i}i}}^{t_{s_{i}i}} z'_{i}v_{\lambda}$$

$$+\sum_{j=2}^{s_{i}} \xi_{ji}t_{ji}z_{(\tilde a_{1i},-m_{1i},0,...,0)}^{t_{1i}}z_{a_{2i}}^{t_{2i}}...z_{a_{ji}+b}z_{a_{ji}}^{t_{ji}-1}z_{a_{j+1}i}^{t_{j+1}i}...z_{a_{s_{i}i}}^{t_{s_{i}i}}z'_{i}v_{\lambda}$$

 \begin{equation}                                                                                                                                     +\sum(\mbox{shorter monomials}). 
\end{equation}

\noindent Since $m_{1i} \geq m_{1i_{0}} > 1$, we see that the first monomial

\begin{eqnarray}
z_{(\tilde a_{1i}+\bar b,-m_{1i}+1,0,...,0)} z_{(\tilde a_{1i},-m_{1i},0,...,0)}^{t_{1i}-1}z_{a_{2i}}^{t_{2i}}...z_{a_{s_{i}i}}^{t_{s_{i}i}}z'_{i} v_{\lambda}      \\ 
\nonumber
\label{eq:m}
\end{eqnarray}
\noindent is still ordered in the ordering defined in Lemma ~\ref{l1}.   
Taking a monomial $M$ from (1), not equal to (2) observe that its smallest term is greater than  
$z_{(\tilde a_{1i}+\bar b, -m_{1i}+1,0,...,0)}$,
 because for $j>1,$ either $m_{ji} > m_{1i}$ or $m_{ji} = m_{1i}$ but $\tilde a_{ji} > \tilde a_{1i}$. 
The monomial $M$ may not be ordered. 
If it is not, then
after ordering $M$ we obtain the ordered permutation of $M$ which is bigger than (2) and some other shorter summands. Thus we can write 
$$yz_{i}v_{\lambda} = \xi_{1i}t_{1i}z_{(\tilde a_{1i}+\bar b,-m_{1i}+1,0,..,0)} z_{(\tilde a_{1i},-m_{1i},0,...,0)}^{t_{1i}-1}z_{a_{2i}}^{t_{2i}}...z_{a_{s_{i}i}}^{t_{s_{i}i}}u'_{i} v_{\lambda} + $$  
\begin{eqnarray}
\sum(\mbox{greater ordered monomials of the same length}) + \sum (\mbox{shorter ordered monomials}),   \\
\nonumber
\label{eq:b}
\end{eqnarray}
\noindent where the first monomial coincides with (2).
 Let $M_{0}$ be the monomial (2)  with $i = i_{0}$.
 So $M_{0}$ is the first monomial in (3) when taking $i=i_{0}.$
Observe that $\xi_{1i_{0}} \neq 0$ as $f(b,a_{1i_{0}}) \neq 1.$ 
It is easy to see that $M_{0}$ does not cancel out in
$$yx = \sum_{i\in I} \lambda_{i}w_{b}z_{i}v_{\lambda}.$$

\noindent Indeed, it follows from (3) that it is enough to check that $M$ does not coincide with (2) when $i \neq i_{0}$. 
But this obviously follows from the choice of $i_{0}$.
 
Observe that $\tilde a_{1i_{0}}+b < 0$ and that $M_{0}$ has degree $\beta =
(\tilde \alpha +\bar b, -L+1,0,...,0)$.
 Notice also that $M_{0}$ does not belong to $\tilde H$. 
This is because $(\tilde a_{1i_{0}}+\bar b,-m_{1i_{0}}+1,0,...,0) \notin R_{r}$ and 
in view of $x \notin \tilde H$, no $a_{ji_{0}} (j = 2,...,s_{i_{0}})$ can be a multiple of some
$a \in \Lambda_{\lambda} \cap R_{r} .$

 Now consider the case when $m_{1i_{0}} = 1.$
Suppose first that there exists $l,$ $1 \leq l \leq s_{i_{0}},$ such that
 $m_{li_{0}} >1.$ Let $l'$ be the least such index. 
By Lemma~\ref{l2}, there exist $N_{1},...,N_{r-1}$ with $N_{j} > k_{j}$ such that 
 $$ f((-N_{1},...,-N_{n-1},1,0,...,0),(\tilde a_{l'i_{0}},-m_{l'i_{0}},0,...,0)) \neq 1.$$
\noindent As before , we take $\bar b = (-N_{1},..., -N_{r-1}),$ $b = (\bar b, 1,0,...,0),$ $y = w_{b}$ and set $\xi_{ji} = \sigma(a_{ji},b)(f(b,a_{ji})-1).$ 
Then $\xi_{l'i_{0}} \neq 0$ and $[w_{b}, z_{a_{ji}}] = \xi_{ji}z_{b+a_{ji}}.$ 
Let $i\in I'$ be an index with similar property: 
$m_{l''i} > 1$ for some  $1 \leq l''\leq s_{i}.$
Suppose that $l''$ is the least such index. 
We have that 
$$yz_{i}v_{\lambda} = \xi_{l''i}t_{l''i}z_{a_{1i}}^{t_{1i}}...z_{a_{l''i}+b}z_{a_{l''i}}^{t_{l''i}-1}...z_{a_{s_{i}i}}^{t_{s_{i}i}}z'_{i}v_{\lambda} + $$
 
$$\sum_{j \neq l''}\xi_{ji}t_{ji}z_{a_{1i}}^{t_{1i}}...z_{a_{ji}+b}z_{a_{ji}}^{t_{ji}-1}...z_{a_{s_{i}i}}^{t_{s_{i}i}}z'_{i} v_{\lambda} + $$
\begin{eqnarray}
\sum (\mbox{shorter monomials}).  \\
\nonumber
\label{eq:c}
\end{eqnarray}
 
\noindent We may still need  to reorder these monomials, but after doing so we will get only shorter summands. 
A routine check shows that among the monomials 
of maximal length in (4) the smallest one is the first summand. 
Let $M_{0}$ be the first monomial of (4) when taking $i = i_{0}$ then ($l''$ becomes
$l'$).
 We claim that $M_{0}$ does not cancel out in $yx.$
We see that it can not happen in $yz_{i}v_{\lambda}$ for every $i \in I'$ with the above property of
existence of $l''.$  
Now, if $i \in I'$ does not have such a property, then the last entries of all $a$'s involved in $z_{i}$ are all $-1$'s 
and $0$'s. Thus when acting on $z_{i}v_{\lambda}$ by $y = w_{(\bar b,1,0,...,0)}$we reduce the number of $-1$'s, while
in $M_{0}$ it is either the same as in $z_{i_{0}}$ (case $m_{l'i_{0}} > 2$), or even greater (case $m_{l'i_{0}} =2$).
Hence, $M_{0}$ does not cancel out in such $yz_{i}v_{\lambda}.$ 
Finally, if $i \notin I'$, then all monomials of $yz_{i}v_{\lambda}$ which appear after reordering are shorter than $M_{0}.$
Thus, $M_{0}$ belongs to the support of $yx$, as $\xi_{l'i_{0}} \neq 0.$ 
Since $M_{0} \notin \tilde H,$ $0 \neq yx \neq \tilde H.$ 
Moreover, the degree of $M_{0}$ is $\beta = (\tilde \alpha +\bar b,-L+1,0,...,0).$

Suppose now that $m_{1i_{0}} = 1$ and $m_{li_{0}} \in \{0,1\}$ for all $1 \leq l \leq s_{i_{0}}.$ 
Let $z_{i_{1}}v_{\lambda}$ be the smallest monomial in $X' \setminus \{ z_{i_{0}}v_{\lambda}\}.$
If again $m_{1i_{1}} =1$ and $m_{li_{1}} \notin \{ 0,1\}$ for all $ 1 \leq l \leq s_{1i},$
then we take the next smallest monomial and continue until we get the first monomial $z_{i_{t}}v_{\lambda}$ with some 
$m_{li_{t}} >1,$ $1 \leq l \leq s_{i_{t}}$.
Let $l'$ be the minimal such index. Let also $l''$ be the maximal index with $m_{l''i_{0}} =1.$ 
Set $I'' = \{i_{0},i_{1},...,i_{t-1} \}$. By Lemma ~\ref{l2} we choose $N_{1},..., N_{r-1}$ with $N_{j} >k_{j}$
such that  $f((-N_{1},...,-N_{r-1},1,0,...,0), a_{l''i_{0}}) \neq 1.$ 
Take $b = (\bar b,1,0,...,0), y = w_{b}.$ Then $\xi_{l''i_{0}} \neq 0$ where $\xi_{ji}$ has the same meaning as before.
For $i \in I''$ we let $l''' $ be
the maximal index with $m_{l'''i} = 1.$
We have $$yz_{i}v_{\lambda} = \xi_{l'''i}t_{l'''i}z_{a_{1i}}^{t_{1i}}...z_{a_{l'''i+b}}z_{a_{l'''i}}^{t_{l'''i}-1}...z_{a_{s_{i}i}}^{t_{s_{i}i}}z'_{i}v_{\lambda} $$
$$+ \sum_{j \neq l'''} \xi_{ji}t_{ji}z_{a_{1i}}^{t_{1i}}...z_{a_{ji}+b}z_{a_{ji}}^{t_{ji}-1}...z_{a_{s_{i}i}}^{t_{s_{i}i}}z'_{i}v_{\lambda}$$
\begin{eqnarray}
\sum(\mbox{shorter monomials}).     \\
\nonumber
\label{eq:d}
\end{eqnarray}
Among the monomials of maximal length in (5) the first summand is the smallest one.
Now let $M_{1}$ be the first monomial of (5) with $i = i_{0}$ (then $l'''$ becomes $l''$).

Since $z_{i_{0}}v_{\lambda} < z_{i}v_{\lambda},$ it follows that $M_{1}$ does not cancel out in $yz_{i}v_{\lambda}$ 
for all $i \in I''.$

If we restrict our attention only to $I' \setminus I'',$ then we are in conditions of previous cases. Let $z_{i'_0}v_{\lambda}$ be the smallest monomial among the monomials $z_iv_{\lambda}$ of maximal length with $i\in I' \setminus I''.$ In view of the previous considerations let $M_{0}$ be the first monomial in (2) or (4) with $i=i'_0$. Now, if $M_{1} < M_{0},$ then $M_{1}$ can not be cancelled out in $yx$ and
it really belongs to the support of $yx$ as $\xi_{l''i_{0}} \neq 0.$\\
Let $M_{0} \leq M_{1}.$ Taking larger $N_1, \ldots, N_{r-1}$ if necessary we can guarantee by Lemma~\ref{l2} that $M_{0} < M_{1}$ and $f(b,a_{l'i'_{0}}) \neq 1$ where, as usually,  $b = (-N_{1},...,-N_{r-1},1,0,...,0).$ Then $\xi_{l'i'_{0}} \neq 0$ and since $M_{0}$ is the smallest one among all monomials of maximal length which appear in $yx,$
it does not cancel out.

%Those considerations show that there is a monomial $M_{0}$ (the first monomial %in (\ref{eq:c}) or (\ref{eq:a})) 
%which is the smallest one among the monomials of maximal length appearing in %$yz_{i}v_{\lambda}$ for all
%$i \in I'\setminus I''.$ Now, if $M_{1} < M_{0},$ then $M_{1}$ can not be %cancelled out in $yx$ and
%it really belongs to the support of $yx$ as $\xi_{l''i_{0}} \neq 0.$\\
%Let $M_{0} < M_{1}.$ By Lemma~\ref{l2}, take $N_{1},...,N_{r-1}$ with $N_{j} > %k_{j}$
%such that $f(b,a_{l'i_{s}}) \neq 1$ where $b = %(-N_{1},...,-N_{n-1},1,0,...,0).$ 
%Then $\xi_{l'i_{s}} \neq 0$ and since $M_{0}$ is the smallest one among all %monomials of maximal length which appear in $yx,$
%it does not cancel out.

\end{proof}
\end{lemma}

\noindent{\bf Proof of Theorem 1.}

\noindent{\it Suppose that $\lambda(c_{1}) \neq 0$. We want to show that  $H(\lambda)/{\tilde H}$ is irreducible.}

 Let $0 \neq W$ be a submodule of $H(\lambda)/{\tilde H}$ and let $0 \neq \bar x = x + \tilde H$ be an element of $W.$
Since $x \notin \tilde H$, we may assume that no monomial in $x$ is in $\tilde H.$ We also may assume that 
 $x$ is homogeneous of degree $\alpha = (\tilde \alpha,-L,0,...,0), \tilde \alpha \in \mathbb Z^{r-1}$  and $L > 0.$ 
 Suppose first that $r = 1.$

 Consider the subalgebra $\mathfrak s$ of $\mathfrak t$ with a $\mathbb C$ - basis 
$\{ c_{1}, w_{(i,0,...,0)}, u_{(j,0,...,0)}| \, i,j \in \mathbb Z \}.$ Then $\mathfrak s$ is a Heisenberg subalgebra of $\mathfrak t$ as in the 
affine setting ( note that $f((i,0,0,...,0),(j,0,0,...,0)) = 1$ for all $i,j \in \mathbb Z$ so that the brackets 
$[w_{(i,0,...,0)},w_{(j,0,...,0)}] = [u_{(i,0,...,0)},u_{(j,0,...,0)}] = 0$ unless $i = -j$  and  in this case
$[w_ {(i,0,...,0)}, w_{(-i,0,...,0)}] = [u_{(i,0,...,0)}, u_{(-i,0,...,0)}] = 2ic_{1}$ and $[u_{(i,0,...,0)}, w_{(j,0,...,0)}] $ $= 0$ for all $i,j \in \mathbb Z.)$
  Since $r=1$, $x \in U(\mathfrak s)v_{\lambda}$ and because $\lambda(c_{1}) \neq 0$, 
we have  from the results
for affine Lie algebras  (see \cite{K}) that $U(\mathfrak s)v_{\lambda}$ is irreducible, so that $v_{\lambda} \in 
U(\mathfrak s)v_{\lambda} \subset U(\mathfrak t)v_{\lambda}.$
  Now, suppose that $r>1$ and that the theorem is true for $r-1.$
  If all  of the monomials in $x$ contain only factors $z_{a}$ for $a \in R_{r}$, 
we proceed exactly as in the proof of [FK, Theorem 1] 
and we obtain $v_{\lambda} \in W.$
   If at least one monomial in $x$ contains a factor $z_{a}$ , with $a = (a_{1},..., a_{r},0,...,0) \notin R_{r}$ and $a_{r} \neq 0$
we proceed by induction on $L$. Applying  Lemma 3 we obtain $0 \neq z \in W$ so that $z \in H(\lambda)_{\beta}$ with 
$\beta = (\tilde \beta,-K,0,....,0)$ where $\tilde \beta \in \mathbb Z^{r-2}$ and $K > 0$ and $z \notin \tilde H.$ The result then 
follows by the induction hypothesis.

\begin{remark}
We believe that Theorem 1 is valid when $\lambda(c_i)\neq 0$ for some $i$.

\end{remark}

%\end{proof}    

\section {Imaginary Verma modules}
      
      In this section we consider certain Verma type modules for the   extended
affine Lie algebra $sl_{2}(\mathbb C_{q}).$ The results are similar to the case
of toroidal Lie algebras [FK].

       Let $$\Delta_{\pm} = \{ \pm \alpha + \Lambda \} \cup \{ a\in \Lambda \, | \, \pm a > 0 \},$$

$$\mathfrak q_{\pm} = \sum_{\beta \in \Delta_{\pm}} \mathfrak q_{\beta} \, \, , \, \mathfrak q_{-}^{re} = \sum _{\beta \in \Delta _{-}\setminus \Lambda} \mathfrak q_{\beta}.$$
           
       We have the following Cartan type decomposition:
           $$\mathfrak q = \mathfrak q_{-} \oplus \mathfrak h \oplus \mathfrak q_{+}.$$ 
We will call a subalgebra $B = \mathfrak h\oplus q_{+}$ a {\it non-standard Borel subalgebra.}

Let $\lambda \in {\mathfrak h}^{*}$. A $\mathfrak q$-submodule 
           $$M(\lambda) = U(\mathfrak q) \otimes _{U(B)}\mathbb C v_{\lambda},$$
\noindent where $\mathfrak q_{+}v_{\lambda} =  0 $ and $h v_{\lambda} = \lambda (h) v_{\lambda},$ is called an 
{\it imaginary Verma module}.

      This is a weight (with respect to $\mathfrak h$) module with $M(\lambda)_{\mu} \neq 0 $ if and only if
$\mu \in \{\lambda - \alpha + \Lambda \} \cup \{ \lambda -a \, | \, a \in \Lambda , a > 0 \}.$  

       Moreover, dim$M(\lambda)_{\lambda} = 1$ and dim$M(\lambda)_{\mu} = \infty$ as long as $M(\lambda)_{\mu} \neq 0$. The module $M(\lambda)$ has a unique irreducible quotient $L(\lambda).$

       Let $\tilde v_{\lambda} = 1 \otimes v_{\lambda}.$ Consider a subspace $\hat M(\lambda) = U(\mathfrak t)\tilde v_{\lambda}.$ Then $\hat M(\lambda)$ is a $\mathfrak t$-submodule of $M(\lambda)$ and  $M(\lambda) \cong H(\lambda).$
       We will show that this submodule plays a crucial role in the structure of $M(\lambda).$ It follows from the PBW theorem that 
          $$M(\lambda) \cong U(\mathfrak q_{-}^{re}) \otimes _{\mathbb C} \hat M(\lambda)$$
\noindent as a vector space.

          For a submodule $N \subset M(\lambda)$ denote $\hat N = N \cap \hat M(\lambda).$        
 \begin{proposition} \label{p3} Let $0 \neq N \subset M(\lambda)$. Then $\hat N \neq 0.$

\begin{proof} Let $0 \neq v \in N$. We can assume that $v$ is a weight vector and that 
$$v = \sum_{i}y_{a_{i1}}...y_{a_{im}}z_{i}\tilde v_{\lambda},$$

\noindent where $a_{ij} \in \Lambda$, $z_{i} \in U(\mathfrak t^{-})$ and all the monomials $y_{a_{i1}}...y_{a_{im}}$
are linearly independent. Note that $m$ is an invariant of $v$ since $v$ is a weight vector. We will denote $\|v\|= m $
and will prove the statement by induction on $m$. The base of induction $m=0$ is trivial. Clearly for appropriate $A \in \Lambda,$
$$N \ni x_{A}v = \sum_{i}\sum_{j=1}^{n}y_{a_{i1}}...[x_{A},y_{a_{ij}}]...y_{a_{im}}z_{i}\tilde v_{\lambda} \neq 0$$
\noindent and $\|x_{A}v\|= m-1$. Hence we can apply our induction hypothesis and conclude that $0 \neq ux_{A}v \in \hat N$
for some $u \in U(\mathfrak q).$

\end{proof}
\end{proposition} 

\begin{corollary}\label{c2}Let $\lambda \in {\mathfrak h}^{*}$. A $\mathfrak q$-module $M(\lambda)$ is irreducible if and
only if $\hat M(\lambda) = H(\lambda)$ is irreducible as a $\mathfrak t$-module.
\end{corollary}

The algebra $U(\mathfrak t)$ has a natural $\Lambda $ -gradation. If $u \in U(\mathfrak t)$ is a
homogeneous element of degree $(a_{1},..., a_{n})$ we denote $|u|_{i} = |a_{i}|.$ For an arbitrary element
 $u \in U(\mathfrak t)$ denote $|u|_{i} = \mbox{max}_{j}|u_{j}|_{i}$ where $u_{j}$ are homogeneous components
 of $u$.

\begin{theorem}\label{t2}Let $N \subset M(\lambda)$. If there exists an 
$i$ such that $\lambda(c_{i}) \neq 0$ then
$N \cong U({\mathfrak q}^{re}_{-})\otimes _{\mathbb C}\hat N$ as vector space.
\end{theorem}

  To prove the theorem we need the following

\begin{lemma}Let $N \subset M(\lambda)$, $\lambda(c_{k}) \neq 0$  and $0 \neq u \tilde v_{\lambda} \in N$ where 
$$u = \sum_{i\in I}z_{A_{i}}u_{i}$$
\noindent with $z \in \{u,w\}$, $A_{i} = (a_{i1},...,a_{in}) < 0$, $A_{i}\notin \Lambda_{\lambda}$, $A_{i} \neq A_{j}$ if $i \neq j$                                   except possibly the case when $z_{A_{i}} = u_{A_{i}}$ and $z_{A_{j}} = w_{A_{j}}$; $u_{i} \in U(\mathfrak t)$ and $|u_{i}|_{j} << |a_{ij}|$,
for all $i \in I$, $j=k$ and $j=n.$ Moreover $|a_{1n}| - 2|u|(n)>> 0$, $|u|(n)=max_k |u_k|_n$. 
 Then $u_{i}\tilde v_{\lambda} \in N$ for all $i \in I.$

\begin{proof}Since $[u_{a}, w_{-a}]= 0$ we can assume without loss of generality that $A_{i} \neq A_{j}$ for all $i \neq j.$ 
Use induction on $I.$  Let $|I|= 1.$ Then $u = z_{A}u'$, $A = (a_{1},..., a_{n}).$ Apply $z_{-A}.$ Since $z_{-A}u'\tilde v_{\lambda}=0$
and $A \notin \Lambda_{\lambda}$ we conclude that $u'\tilde v_{\lambda} \in N.$

Suppose now that $|I| > 1, |a_{1n}| \leq |a_{2n}| \leq ...$ and $(a_{11},...,a_{1n-1}) > (a_{i1},..., a_{in-1})$ if $a_{in} = a_{1n} .$
We can also assume that $f(A_{1}, A_{i}) \neq 1$ for some $i \neq 1$ since otherwise application of $z_{-A_{1}}$ completes the 
proof. Consider such $i$. Assume that $k\neq n$ and consider 
$$B_{D,M} = (-a_{11},...,-a_{1k-1}, -a_{1k}-D,-a_{1k+1},..., -a_{1n-1}, 
-a_{1n}-M),$$
where $D>>0$ and $|u|(n)<< M << |a_{1n}|- |u|(n)$.
\noindent Suppose that $f(B_{D,M},A_{i}) = 1$, for sufficiently large number of pairs $(D,M)$. Then 

$$(q_{k1}^{a_{i1}}q_{k2}^{a_{i2}}...q_{kk-1}^{a_{ik-1}}q_{kk+1}^{a_{ik+1}}...q_{kn}^{a_{in}})^D(q_{n1}^{a_{i1}}...q_{nn-1}^{a_{in-1}})^M$$
\noindent is a constant for all those pairs $(D,M)$. Therefore 
$$ q_{k1}^{a_{i1}}q_{k2}^{a_{i2}}...q_{kn}^{a_{in}}= q_{n1}^{a_{i1}}...q_{nn-1}^{a_{in-1}} =1$$
\noindent implying that $f(-A_{1}, A_{i}) =1$ which is a contradiction. Thus we can assume that
  $$f(B_{D,M},A_{i}) \neq 1$$ for a large
number of pairs $(D,M)$ and some $i\neq 1.$ We can also choose $D$ and $M$ in such a way that 
$\tilde A_{i} = A_{i}+B_{D,M} \notin \Lambda_{\lambda}$ for all $i \in I.$ Fix $D$ and $M$ that
satisfy conditions above and  apply $u_{B_{D,M}}$:

$$v= u_{B_{D,M}}u\tilde v_{\lambda} = [u_{B_{D,M}},z_{A_{1}}]u_{1}\tilde v_{\lambda} + 
\sum_{i \neq 1}[u_{B_{D,M}}, z_{A_{i}}]u_{i}\tilde v_{\lambda}
\in N.$$

\noindent If $f(B_{D,M},A_{1})= 1 $ then the first term in the sum above is zero 
and we can complete the proof by induction on $|I|.$ Assume now that $f(B_{D,M}, A_{1}) \neq 1$.
 Denote by $\tilde z_{i} = [u_{B_{D,M}}, z_{A_{i}}], \tilde A_{i} = 
(\tilde a_{i1},...,\tilde a_{in}).$ If $f(-\tilde A_{1}, \tilde A_{i}) = 1$
for all $i \neq 1$ then $z_{-\tilde A_{1}}v = [z_{-\tilde A_{1}},\tilde z_{1}]u_{1}\tilde v_{\lambda}$, 
with $z$ here the same as in $\tilde z_{1}$,  
and therefore $u_{1}\tilde v_{\lambda} \in N.$ Thus we assume that  
 $f(-\tilde A_{1}, \tilde A_{i}) \neq  1$ for some $i \neq 1.$ Consider $C_{K,L} = (C_{1},...,C_{n})$ 
where $C_{k} = D-K > > 0,$
$C_{n} = M-L > > 0,$ $C_{j} =0$ for $j \neq k,n.$ The same argument as above shows that there exists a sufficient number of pairs
$(K,L)$ such that $f(C_{K,L}, \tilde A_{i}) \neq 1$ for some $i\neq 1,$ and for all such indices 
$i$, $C_{K,L} + \tilde A_{i} \notin \Lambda_{\lambda}.$ First note that among all 
 pairs $(K,L)$ there are some for which $f(C_{K,L}, \tilde A_{1}) =1.$ Indeed, 
$f(C_{K,L}, \tilde A_{1}) =1$ as soon as $KM =LD.$ Without loss of generality we can assume that $M$ and $D$ 
are both divisible by $2^r$ for some integer $r>>0$.
Set $K = \frac {1}{2^r}D$, $L=\frac {1}{2^r}M$ and assume that 
$f(C_{\frac {1}{2^r}D,\frac {1}{2^r}M},\tilde A_{i}) =1.$ for sufficiently 
 many integers $r$. 
 Hence

$$q_{k1}^{\tilde a_{i1}(D-\frac{1}{2^{r}}D)}...q_{kn}^{\tilde a_{in}(D-\frac{1}{2^{r}}D)}q_{n1}^{\tilde a_{i1}(M-\frac{1}{2^{r}}M)}...q_{n n-1}^{\tilde a_{in-1}(M-\frac{1}{2^{r}}M)} =1$$
\noindent and $$f(-\tilde A_{1}, \tilde A_{i}) = 
(q_{k1}^{\tilde a_{i1}}...q_{kn}^{\tilde a_{in}})^{D}(q_{n1}^{\tilde a_{i1}}...q_{nn-1}^{\tilde
a_{in-1}})^{M} = 1.$$
\noindent We conclude that there exist infinitely many pairs $(K,L)$ satisfying $KM=LD$ and $f(C_{K,L}, \tilde A_{i}) \neq 1$ as long
as $f(-\tilde A_{1}, \tilde A_{i}) \neq 1.$ Finally, if $C_{\frac {1}{2^r}D,\frac {1}{2^r}M} +\tilde A_{i}  \in \Lambda_{\lambda}$
then
$$0 = ((D-\frac{1}{2^{r}})+\tilde a_{ik})\lambda(c_{k}) + (M-\frac{1}{2^r}M)\lambda(c_{n}) + \sum_{j\neq k,n} \tilde a_{ij}\lambda (c_{j}) =$$
$$(1-\frac{1}{2^r})(D\lambda (c_{k})+M\lambda (c_{n})) + \sum_{j}\tilde a_{ij} \lambda(c_{j}),$$
\noindent which  is impossible since $\tilde A_{i} \notin \Lambda_{\lambda}$ for all $i\in I.$ 
Hence, we can choose $0< < K < < D$ and
$0 < < L < < M$ such that $f(C_{K,L}, \tilde A_{1}) = 1$ and $f(C_{K,L}, \tilde A_{i}) \neq 1$ for at least one $i\in I$, and in this case $C_{K,L} + \tilde A_{i} \notin \Lambda_{\lambda}.$
Apply $u_{C_{K,L}}$:
$$0 \neq  u_{C_{K,L}}v = \sum_{i\neq 1}[u_{C_{K,L}}, \tilde z_{i}]u_{i}\tilde v_{\lambda} \in N.$$
 We complete the proof by induction on $|I|.$ The same arguments work in the case $k=n$. 
The lemma is proved.
\end{proof}
\end{lemma}

The Theorem 2 is an easy corollary  of Lemma 4. The proof follows the general lines of the proof of Theorem 2 in
\cite{FK}. We omit the details here.

\begin{corollary}
Let $\lambda \in {\mathfrak h}^*$ and $\lambda(c_i)\neq 0$ for some $i$. Then   $L(\lambda) \cong 
U({\mathfrak q}^{re}_{-})
\otimes_{\mathbb C} \bar{H}(\lambda)$, where $\bar{H}(\lambda)$ is the maximal submodule of $H(\lambda)$.

\end{corollary}

\begin{corollary}
Let $\lambda \in {\mathfrak h}^*$ and $\lambda(c_1)\neq 0$. 
Then $L(\lambda) \cong U({\mathfrak q}^{re}_{-})
\otimes_{\mathbb C} (H(\lambda)/{\tilde H})$.
\end{corollary}
\begin{proof}
Follows from Theorem 1 and Theorem 2. 
\end{proof}

\begin{center}
{\bf\large Acknowledgements}
\end{center}
This work was partially done during the stay of the third author at the University of S\~ao Paulo as a visiting 
professor and was completed during his visit to the Fields Institute. 
The financial support of  Fapesp of Brazil, Proc. 97/05415-0,  and the Fields Institute is gratefully acknowledged.
 The third author is grateful to B.Allison for 
many useful discussions on the subject. The first author was partially supported by CNPq of Brazil, Proc. 301115/95-8.

\vspace{1cm}

\noindent

\noindent Instituto de Matem\'atica e Estatistica,\\
Universidade de S\~ao Paulo\\
Caixa Postal 66281- CEP 05315-970\\
S\~ao Paulo, Brazil\\
e-mail: futorny@ime.usp.br

\end{document}